\documentclass[a4paper,11pt]{amsart}

\usepackage{amsmath}
\usepackage{amsthm}
\usepackage[all]{xy}
\usepackage{amssymb}
\usepackage{mathrsfs}
\usepackage{enumerate}
\usepackage[usenames,dvipsnames]{color}

\usepackage{fullpage}
\usepackage{amsmath,amsfonts,amssymb}
\makeatletter
\newcommand{\thorn}{{\fontencoding{T1}\selectfont\th}}
\newcommand{\@indepsymbol}[2]{#1\setbox0=\hbox{$#1x$}\kern\wd0\hbox to 0pt{\hss$#1\mid$\hss}\lower.9\ht0\hbox to 0pt{\hss$#1\smile$\hss}\kern\wd0}
\newcommand{\@nindepsymbol}[2]{#1\setbox0=\hbox{$#1x$}\kern\wd0\hbox to 0pt{\mathchardef
	\nn=12854\hss$#1\nn$\kern1.4\wd0\hss}\hbox to
	0pt{\hss$#1\mid$\hss}\lower.9\ht0 \hbox to
	0pt{\hss$#1\smile$\hss}\kern\wd0}
\newcommand{\ind}[1][]{\mathop{\mathpalette\@indepsymbol{}^{\!\!\!\!\rlap{$\scriptstyle\textnormal{#1}$}\,\,\,\,}}}
\newcommand{\nind}[1][]{\mathop{\mathpalette\@nindepsymbol{}^{\!\!\!\rlap{$\scriptstyle\textnormal{#1}$}\,\,\,}}}
\newcommand{\@Ind}[1][]{\mathpalette\@indepsymbol{}^{\!\!\!\!\mbox{$\scriptstyle\textnormal{#1}$}}}
\newcommand{\Ind}[1][]{\@Ind[\ \,]}
\newcommand{\newind}[4]{
	\newcommand{#1}{{\!\@Ind[#4]}}
	\newcommand{#2}{\ind[#4]}
	\newcommand{#3}{\nind[#4]}
}
\newind{\thInd}{\thind}{\nthind}{\thorn}
\renewcommand{\thInd}{\text{$\@Ind[\thorn]$\;}}
\makeatother

\newtheorem{thm}{Theorem}[section]
\newtheorem{cor}[thm]{Corollary}
\newtheorem{conj}[thm]{Conjecture}
\newtheorem{lem}[thm]{Lemma}
\newtheorem{prop}[thm]{Proposition}

\newtheorem{fact}[thm]{Fact}
\theoremstyle{definition}
\newtheorem{defn}[thm]{Definition}
\theoremstyle{remark}

\theoremstyle{remark}

\theoremstyle{remark}

\title{Points on a curve with a power on a curve}
\date{\today}
\author{Gareth J. Boxall}
\thanks{This work is based on the research supported in part by the National Research Foundation of South Africa (Grant Number 96234).}
\email{gboxall@sun.ac.za, garethjb@gmail.com}
\address{Department of Mathematical Sciences (Mathematics Division), Stellenbosch University, Private Bag X1, Matieland 7602, South Africa}

\begin{document}
\maketitle

\begin{abstract}
Let $C_1,C_2\subseteq\mathbb{G}_m^N(\mathbb{C})$ be irreducible closed algebraic curves, with $N\geq 3$. Suppose $C_1$ is not contained in an algebraic subgroup of $\mathbb{G}_m^N(\mathbb{C})$ of dimension $1$ and $C_1\cup C_2$ is not contained in an algebraic subgroup of $\mathbb{G}_m^N(\mathbb{C})$ of dimension $2$. It is a conjecture that at most finitely many points $x\in C_1$ have the property that there is a positive integer $n$ such that $x^n\in C_2$ and $[n]C_1\nsubseteq C_2$, where $[n]C_1=\{x^n:x\in C_1\}$. We prove this in the case where at least one of the two curves is not defined over $\overline{\mathbb{Q}}$. 
\end{abstract}

\section{Introduction}

We consider the following conjecture.

\begin{conj}\label{stronger}
Let $C_1,C_2\subseteq\mathbb{G}_m^N(\mathbb{C})$ be irreducible closed algebraic curves, with $N\geq 3$. Suppose there does not exist an algebraic subgroup $G\subseteq \mathbb{G}_m^N(\mathbb{C})$ of dimension one such that $C_1\subseteq G$ and there does not exist an algebraic subgroup $H\subseteq\mathbb{G}_m^N(\mathbb{C})$ of dimension two such that $C_1\cup C_2\subseteq H$. Let $\mathcal{N}=\{n\in\mathbb{N}:[n]C_1\subseteq C_2\}$. Then $\bigcup\limits_{n\in \mathbb{N}\setminus\mathcal{N}}\{x\in C_1:x^n\in C_2\}$ is finite. 
\end{conj}

Here and throughout, for subvarieties of $\mathbb{G}_m^N(F)$, ``irreducible" means ``irreducible over $F$". Addressing a question of Aaron Levin, Zannier showed in \cite{Z} how to deduce Conjecture \ref{stronger} from Zilber's conjecture on unlikely intersections with tori (CIT). We prove the following.

\begin{thm}\label{newtheorem}
Conjecture \ref{stronger} is true under the additional assumption that least one of the two curves is not defined over $\overline{\mathbb{Q}}$.
\end{thm}

This is not dependent on the case where $C_1$ and $C_2$ are defined over $\overline{\mathbb{Q}}$, which remains open. To prove Theorem \ref{newtheorem} we use work of Bays, Kirby and Wilkie \cite{BKW} which established an analogue of Schanuel's conjecture for the operation of raising to an exponentially transcendental power. Their reasoning relied on separate work of Kirby \cite{K} on partial exponential fields which, in turn, made use of Ax's famous paper \cite{A} on functional analogues of Schanuel's conjecture. The main result of \cite{BKW} is the following.

\begin{thm}\label{BKW}
Let $F$ be an exponential field and $N$ a positive integer. Let $a\in F^N$, $\lambda\in F$ and $K\subseteq F$ be such that $K$ is a subfield, $a$ is $\mathbb{Q}$-linearly independent and $\lambda$ is exponentially transcendental over $K$. The field generated by $(a,\exp(a))$ has transcendence degree at least $N$ over $K(\lambda)$. 
\end{thm}

An exponential field is a field of characteristic zero equipped with a homomorphism $\exp$ from the additive group to the multiplicative group. Exponential transcendence is a natural analogue, involving $\exp$, of the usual notion of transcendence. For our argument, we require a slight generalisation of Theorem \ref{BKW} to the setting of partial exponential fields, where $\exp$ might not be defined everywhere. Given the results in \cite{K}, it is easy to adapt the reasoning in \cite{BKW} to obtain this. We apply it in a structure formed by taking an elementary extension of the complex exponential field, restricting $\exp$ to the infinitesimal elements and then forming an ultrapower of the resulting structure. Our exponentially transcendental exponents are non-standard positive integers. 

Aspects of the argument work for higher dimensional subvarieties and we obtain the following.

\begin{thm}\label{higher dim}
Let $V_1, V_2\subseteq\mathbb{G}_m^N(\mathbb{C})$ be irreducible closed algebraic subvarieties, with $\dim(V_1)\geq 1$ and $\dim(V_1)+\dim(V_2)< N$. Suppose that, for any $\mathbb{Q}$-linearly independent tuples of integers $$(n_{1,1},...,n_{1,N}),...,(n_{k,1},...,n_{k,N}),$$ with $k=\dim(V_1)+1$, and any generic point $(x_1,...,x_N)$ of $V_1$, the monomials $$x_1^{n_{1,1}}\cdot...\cdot x_N^{n_{1,N}},...,x_1^{n_{k,1}}\cdot...\cdot x_N^{n_{k,N}}$$ are algebraically independent.  Let $\mathcal{N}=\{n\in\mathbb{N}:[n]V_1\subseteq V_2\}$. Then $\bigcup\limits_{n\in\mathbb{N}\setminus\mathcal{N}}\{x\in V_1:x^n\in V_2\}$ is not Zariski-dense in $V_1$. 
\end{thm}

A similar algebraic independence condition on multiplicative combinations of coordinates was used by Habegger in \cite{HBH} to isolate a class of subvarieties for which he could prove a version of CIT, though there $k$ was equal to the dimension and the subvariety was assumed to be definable over $\overline{\mathbb{Q}}$. Masser considered subvarieties satisfying the assumption for $k=\dim(V)+1$ in positive characteristic in \cite{Mass}. 

Interest in CIT motivated an investigation by Bays and Habegger into Levin's question in \cite{BH}. They considered the case of a single curve $C=C_1=C_2$ and obtained the following where, with $x=(x_1,...,x_N)$, $$\mathbb{T}=\{x\in\mathbb{G}_m^N(\mathbb{C}):\left|x_1\right|=...=\left|x_N\right|=1\}.$$

\begin{thm}\label{BH}
Let $C\subseteq\mathbb{G}_m^N(\mathbb{C})$ be an irreducible closed algebraic curve defined over $\overline{\mathbb{Q}}$, with $N\geq 3$. Suppose $C$ is not contained in a proper algebraic subgroup of $\mathbb{G}_m^N(\mathbb{C})$. \begin{enumerate} \item If $\sigma(C)\cap \mathbb{T}$ is finite, for some automorphism $\sigma$ of $\mathbb{C}$, then $\bigcup\limits_{n\in \mathbb{N}\setminus\{1\}}\{x\in C:x^n\in C\}$ is finite. \item If $S$ is a finite set of sufficiently large prime numbers then $\bigcup\limits_{n\in \mathbb{N}\setminus\{1\}}\{x\in C:x^n\in C\}$ contains at most finitely many $S$-integral points. \end{enumerate}
\end{thm}

As is noted in \cite{BH}, it follows from known results that $\{n\in\mathbb{N}:[n]C\subseteq C\}=\{1\}$ for such $C$. So Theorem \ref{BH} gives two special cases of a statement that generalises to the two curve setting as follows.

\begin{conj}\label{main}
Let $C_1,C_2\subseteq\mathbb{G}_m^N(\mathbb{C})$ be irreducible closed algebraic curves, with $N\geq 3$. Assume there is no proper algebraic subgroup $H\subseteq \mathbb{G}_m^N(\mathbb{C})$ such that $C_1\subseteq H$. Let $\mathcal{N}=\{n\in\mathbb{N}:[n]C_1\subseteq C_2\}$. Then $\bigcup\limits_{n\in \mathbb{N}\setminus\mathcal{N}}\{x\in C_1:x^n\in C_2\}$ is finite. 
\end{conj}

Conjecture \ref{stronger} follows easily from Conjecture \ref{main}, given known results. Other important work towards Conjecture \ref{main}, for curves defined over $\overline{\mathbb{Q}}$, was done by Amoroso, Masser and Zannier in \cite{AMZ}. They showed that the points in $\bigcup\limits_{n\in \mathbb{N}\setminus\mathcal{N}}\{x\in C_1:x^n\in C_2\}$ will have bounded height. Using this it is possible to extend Theorem \ref{BH} to the two curve setting. I shall include that in a subsequent paper (or an expansion of this one).

This paper is arranged as follows. In Section 2 we recall some known results concerning unlikely intersections in $\mathbb{G}_m^N(\mathbb{C})$, including the dimension one case of CIT and a theorem of Liardet. We also present a deduction of Conjecture \ref{stronger} from Conjecture \ref{main}. In Section 3 we recall some facts about partial exponential fields, from \cite{K}, and slightly modify the proof of the main result in \cite{BKW} to obtain the generalisation that we shall want. In Section 4 we use this to prove the case of Conjecture \ref{main} in which $C_1$ or $C_2$ is not defined over $\overline{\mathbb{Q}}$. By reasoning presented in Section 2, this implies Theorem \ref{newtheorem}. We also prove Theorem \ref{higher dim}. 

I am grateful to Martin Bays, Philipp Habegger, Andrew Harrison-Migochi, Gareth Jones, Sophie Marques, Amador Martin-Pizarro, Harry Schmidt and David Smith for some helpful comments and conversations while I was working on this material. I especially thank Gareth for suggesting that some of the reasoning in the proof of Theorem \ref{newtheorem} would extend to higher dimensional subvarieties, which led to Theorem \ref{higher dim}. 

\section{CIT for curves and a theorem of Liardet}

In this section we discuss some preliminaries. We begin with basic information about algebraic subgroups of $\mathbb{G}_m^N(F)$, for an arbitrary field $F$ of characteristic zero, which can be found in Chapter 3 of \cite{BG}. 

\begin{fact}\label{alggroup}
Let $N$ be an integer greater than or equal to two. Let $F$ be a field of characteristic zero.

\begin{enumerate}

\item A subset $H\subseteq \mathbb{G}_m^N(F)$ is an algebraic subgroup if and only if it is the solution set of a finite system of equations of the form $X_1^{j_{1,1}}\cdot ... \cdot X_N^{j_{1,N}}=1,\ ... ,\ X_1^{j_{s,1}}\cdot ...\cdot X_N^{j_{s,N}}=1$, where $j_{i,q}\in \mathbb{Z}$ for all $i\in\{1,...,s\}$ and $q\in \{1,...,N\}$. 

\item If $H$ is an algebraic subgroup of $\mathbb{G}_m^N(F)$ and the exponents in the equations are as above then $\dim(H)=N-D$, where $D$ is the dimension of the $\mathbb{Q}$-vector space generated by $(j_{1,1},...,j_{1,N}), ... , (j_{s,1},...,j_{s,N})$. 

\end{enumerate}

\end{fact}

The following result is the dimension one case of CIT. It was proved by Maurin \cite{M} for curves defined over $\overline{\mathbb{Q}}$. Bombieri, Masser and Zannier \cite{BMZ} extended it to curves defined over $\mathbb{C}$. 

\begin{thm}\label{MT}
Let $C\subseteq \mathbb{G}_m^N(\mathbb{C})$ be an irreducible closed algebraic curve. Assume $N\geq 2$ and that $C$ is not contained in a proper algebraic subgroup of $\mathbb{G}_m^N(\mathbb{C})$. At most finitely many points $x\in C$ have the property that $x$ is contained in an algebraic subgroup of $\mathbb{G}_m^N(\mathbb{C})$ of dimension $N-2$. 
\end{thm}

The following result was proved by Liardet \cite{Li}, verifying a conjecture of Lang and building on work of Lang and others. It is known to be a consequence of Theorem \ref{MT}, but was established earlier. 

\begin{thm}\label{Liardet1}
Let $C\subseteq \mathbb{G}_m^N(\mathbb{C})$ be an irreducible closed algebraic curve, with $N\geq 2$. Let $\Gamma$ be a finite rank subgroup of $\mathbb{G}_m^N(\mathbb{C})$. If $C\cap \Gamma$ is infinite then $C$ is a translate of an algebraic subgroup of $\mathbb{G}_m^N(\mathbb{C})$. 
\end{thm}

In some important special cases, the conclusion can be strengthened. The following is a well-known consequence of Theorem \ref{Liardet1}. We include a proof for convenience. In fact, the second and third parts of the result were known before Theorem \ref{Liardet1} (see \cite{Lang1} and \cite{Lang2}). 

\begin{cor}\label{Liardet}
Let $C\subseteq \mathbb{G}_m^N(\mathbb{C})$ be an irreducible closed algebraic curve, with $N\geq 2$. Assume at least one of the following three conditions is satisfied.

\begin{enumerate}

\item There is a non-torsion $x\in \mathbb{G}_m^N(\mathbb{C})$ such that, for infinitely many $n\in \mathbb{N}$, $x^n\in C$.

\item There is a non-torsion $y\in \mathbb{G}_m^N(\mathbb{C})$ such that, for infinitely many $n\in \mathbb{N}$, there exists $x\in C$ such that $x^n=y$.

\item $C$ contains infinitely many torsion points.

\end{enumerate} Then $C$ is contained in an algebraic subgroup of $\mathbb{G}_m^N(\mathbb{C})$ of dimension $1$. 
\end{cor}

\proof Assume (1), (2) or (3). By Theorem \ref{Liardet1}, $C$ is a translate of an algebraic subgroup of $\mathbb{G}_m^N(\mathbb{C})$ of dimension $1$. So $C$ is contained in $N-1$ hypersurfaces, defined by equations of the form $X_1^{i_1}\cdot...\cdot X_N^{i_N}=\alpha$, and it is sufficient to show that, in each case, $\alpha$ is a root of unity. If (1) holds then we obtain $\beta\in \mathbb{C}\setminus\{0\}$ and distinct $n,m\in \mathbb{N}$ such that $\beta^n=\alpha=\beta^m$. If (2) holds then we obtain $\beta\in \mathbb{C}\setminus\{0\}$ and distinct $n,m\in\mathbb{N}$ such that $\alpha^n=\beta=\alpha^m$. If (3) holds then $\alpha$ is a multiplicative combination of roots of unity. In each case, $\alpha$ is a root of unity. \endproof

The following result and its proof are also known (see Section 2.4 of \cite{AMZ}). 

\begin{cor}\label{N is finite}
Let $C_1,C_2\subseteq\mathbb{G}_m^N(\mathbb{C})$ be irreducible closed algebraic curves. Assume $N\geq 2$ and that there is no algebraic subgroup $H\subseteq \mathbb{G}_m^N(\mathbb{C})$ of dimension $1$ such that $C_1\subseteq H$. Let $\mathcal{N}=\{n\in\mathbb{N}:[n]C_1\subseteq C_2\}$. Then $\mathcal{N}$ is finite.
\end{cor}

\proof Suppose $\mathcal{N}$ is infinite. Let $x\in C_1$ be generic over a countable subfield $K\subseteq\mathbb{C}$ over which both $C_1$ and $C_2$ are defined. Since $x$ is not a torsion point, it follows by Corollary \ref{Liardet} that $C_2$ is contained in an algebraic subgroup of $\mathbb{G}_m^N(\mathbb{C})$ of dimension $1$. Then the multiplicative dimension of $x$ is at most $1$. This is a contradiction. \endproof

We can deal quickly with one special case of Conjecture \ref{main} and it is convenient to do so. In the proof of the following, $\pi_{1,2}:\mathbb{G}_m^N(\mathbb{C})\rightarrow\mathbb{G}_m^2(\mathbb{C})$ is the projection to the first two coordinates. 

\begin{prop}\label{C not a coset}
Let $C_1,C_2\subseteq\mathbb{G}_m^N(\mathbb{C})$ be irreducible closed algebraic curves, with $N\geq 3$. Assume $C_1$ is not contained in a proper algebraic subgroup of $\mathbb{G}_m^N(\mathbb{C})$. Let $\mathcal{N}=\{n\in\mathbb{N}:[n]C_1\subseteq C_2\}$. If $C_1$ or $C_2$ is contained in a translate of an algebraic subgroup of $\mathbb{G}_m^N(\mathbb{C})$ of dimension $N-2$ then $\bigcup\limits_{n\in \mathbb{N}\setminus\mathcal{N}}\{x\in C_1:x^n\in C_2\}$ is finite. 
\end{prop}

\proof Suppose $\bigcup\limits_{n\in \mathbb{N}\setminus\mathcal{N}}\{x\in C_1:x^n\in C_2\}$ is infinite and $C_1$ is contained in a translate of an algebraic subgroup of $\mathbb{G}_m^N(\mathbb{C})$ of dimension $N-2$. We may assume there exist $\alpha,\beta\in \mathbb{C}\setminus\{0\}$ such that $\pi_{1,2}(x)=(\alpha,\beta)$ for all $x\in C_1$. Then $(\alpha,\beta)$ is not a torsion point. By Corollary \ref{Liardet}, $\pi_{1,2}(C_2)$ is contained in an algebraic subgroup of $\mathbb{G}_m^2(\mathbb{C})$ of dimension $1$. Then $\alpha$ and $\beta$ are multiplicatively dependent, contradicting the assumption that $C_1$ is not contained in a proper algebraic subgroup of $\mathbb{G}_m^N(\mathbb{C})$. 

Suppose $\bigcup\limits_{n\in \mathbb{N}\setminus\mathcal{N}}\{x\in C_1:x^n\in C_2\}$ is infinite and $C_2$ is contained in a translate of an algebraic subgroup of $\mathbb{G}_m^N(\mathbb{C})$ of dimension $N-2$. We may assume there exist $\alpha,\beta\in \mathbb{C}\setminus\{0\}$ such that $\pi_{1,2}(y)=(\alpha,\beta)$ for all $y\in C_2$. Given what we have already proved, we may assume $\pi_{1,2}(C_1)$ is infinite. Then $\bigcup\limits_{n\in \mathbb{N}\setminus\mathcal{N}}\{x\in \pi_{1,2}(C_1):x^n\in \pi_{1,2}(C_2)\}$ is infinite. By Corollary \ref{Liardet}, $\pi_{1,2}(C_1)$ is contained in a proper algebraic subgroup of $\mathbb{G}_m^2(\mathbb{C})$. Therefore $C_1$ is contained in a proper algebraic subgroup of $\mathbb{G}_m^N(\mathbb{C})$, which is a contradiction. \endproof

We conclude this section with details of the deduction of Conjecture \ref{stronger} from the $N=3$ case of Conjecture \ref{main}. 

\begin{prop}\label{implies}
Suppose it is the case that if $C'_1,C'_2\subseteq\mathbb{G}_m^3(\mathbb{C})$ are irreducible closed algebraic curves and $C'_1$ is not contained in a proper algebraic subgroup of $\mathbb{G}_m^3(\mathbb{C})$ then, with $\mathcal{N}'=\{n\in\mathbb{N}:[n]C'_1\subseteq C'_2\}$, the set $\bigcup\limits_{n\in\mathbb{N}\setminus\mathcal{N}'}\{x\in C'_1:x^n\in C'_2\}$ is finite. Then it is also the case that if $C_1,C_2\subseteq\mathbb{G}_m^N(\mathbb{C})$ are irreducible closed algebraic curves, with $N\geq 3$, and there does not exist an algebraic subgroup $G\subseteq\mathbb{G}_m^N(\mathbb{C})$ of dimension $1$ such that $C_1\subseteq G$ and there does not exist an algebraic subgroup $H\subseteq \mathbb{G}_m^N(\mathbb{C})$ of dimension $2$ such that $C_1\cup C_2\subseteq H$ then, with $\mathcal{N}=\{n\in\mathbb{N}:[n]C_1\subseteq C_2\}$, the set $\bigcup\limits_{n\in\mathbb{N}\setminus\mathcal{N}}\{x\in C_1:x^n\in C_2\}$ is finite. 
\end{prop}

\proof Let $N\geq 3$ and let $C_1,C_2\subseteq\mathbb{G}_m^N(\mathbb{C})$ be irreducible closed algebraic curves such that $C_1$ is not contained in an algebraic subgroup of $\mathbb{G}_m^N(\mathbb{C})$ of dimension $1$ and $C_1\cup C_2$ is not contained in an algebraic subgroup of $\mathbb{G}_m^N(\mathbb{C})$ of dimension $2$. Let $\mathcal{N}=\{n\in\mathbb{N}:[n]C_1\subseteq C_2\}$. Suppose $\bigcup\limits_{n\in\mathbb{N}\setminus\mathcal{N}}\{x\in C_1:x^n\in C_2\}$ is infinite. We consider two cases and obtain a contradiction in each.

Case 1: Suppose $C_1$ is not contained in an algebraic subgroup of $\mathbb{G}_m^N(\mathbb{C})$ of dimension $2$. Then there is a projection $\pi:\mathbb{G}_m^N(\mathbb{C})\rightarrow\mathbb{G}_m^3(\mathbb{C})$ such that $\pi(C_1)$ is not contained in a proper algebraic subgroup of $\mathbb{G}_m^3(\mathbb{C})$. Let $C'_1$ and $C'_2$ be the Zariski closures of $\pi(C_1)$ and $\pi(C_2)$, respectively, in $\mathbb{G}_m^3(\mathbb{C})$. These are both irreducible. If either $C'_1$ or $C'_2$ is finite then, using Corollary \ref{Liardet} when one of them is infinite, we obtain the contradictory conclusion that $C'_1$ is contained in a proper algebraic subgroup of $\mathbb{G}_m^3(\mathbb{C})$. So $C'_1$ and $C'_2$ are curves. Let $\mathcal{N}'=\{n\in\mathbb{N}:[n]C'_1\subseteq C'_2\}$. By assumption, $\bigcup\limits_{n\in \mathbb{N}\setminus\mathcal{N}'}\{x\in C'_1:x^n\in C'_2\}$ is finite. Since $C_1$ is irreducible and $C_1'$ is infinite, it then follows that $\bigcup\limits_{n\in \mathbb{N}\setminus\mathcal{N}'}\{x\in C_1:x^n\in C_2\}$ is finite. By Corollary \ref{N is finite}, $\mathcal{N}'$ is finite. Therefore $\bigcup\limits_{n\in \mathbb{N}\setminus\mathcal{N}}\{x\in C_1:x^n\in C_2\}$ is finite.

Case 2: Suppose there is an algebraic subgroup $H\subseteq\mathbb{G}_m^N(\mathbb{C})$ of dimension $2$ such that $C_1\subseteq H$. Given that $\bigcup\limits_{n\in \mathbb{N}\setminus\mathcal{N}}\{x\in C_1:x^n\in C_2\}$ is infinite, there are two subcases to consider.

Subcase 2(a): Assume infinitely many $y\in C_2$ have the property that there exist $x\in C_1$ and $n\in\mathbb{N}\setminus\mathcal{N}$ such that $x^n=y$. Then $C_2\subseteq H$.

Subcase 2(b): Assume only finitely many $y\in C_2$ have this property. Then there exists $y\in C_2$ such that, for infinitely many $x\in C_1$, there is some $n\in\mathbb{N}\setminus\mathcal{N}$ such that $x^n=y$. By Corollary \ref{Liardet}, $C_1$ is contained in an algebraic subgroup of $\mathbb{G}_m^N(\mathbb{C})$ of dimension $1$.

In either subcase we get a contradiction. \endproof

We shall want a relative version of this, given that we are interested in curves that are not both defined over $\overline{\mathbb{Q}}$.

\begin{prop}\label{relative implies}
The case of Conjecture \ref{main} in which $N=3$ and at least one of the two curves is not defined over $\overline{\mathbb{Q}}$ implies the case of Conjecture \ref{stronger} in which $N\geq 3$ and at least one of the two curves is not defined over $\overline{\mathbb{Q}}$. 
\end{prop}

\proof We can use almost exactly the same proof as for Proposition \ref{implies}. We just have to check that if $C_1$ or $C_2$ is not defined over $\overline{\mathbb{Q}}$ then the projection $\pi:\mathbb{G}_m^N(\mathbb{C})\rightarrow\mathbb{G}_m^3(\mathbb{C})$ in Case 1 can be chosen so that $C'_1$ or $C'_2$ is not defined over $\overline{\mathbb{Q}}$. Suppose this is not the case. Let $\pi':\mathbb{G}_m^N(\mathbb{C})\rightarrow\mathbb{G}_m^k(\mathbb{C})$ be a projection such that $\pi'(C_1)$ is not contained in a proper algebraic subgroup of $\mathbb{G}_m^k(\mathbb{C})$. Assume $k$ is the maximum integer for which this is possible. Then $k\geq 3$. By Corollary \ref{Liardet}, $\bigcup\limits_{n\in\mathbb{N}\setminus\mathcal{N}}\{y\in C_2:\text{there exists }x\in C_1\text{ such that }x^n=y\}$ is infinite and so, for each $y\in C_2$, the multiplicative dimension of $y$ over $\pi'(y)$ is zero. By our assumption, the Zariski closures of $\pi'(C_1)$ and $\pi'(C_2)$ in $\mathbb{G}_m^k(\mathbb{C})$ are both defined over $\overline{\mathbb{Q}}$. It follows that both $C_1$ and $C_2$ are defined over $\overline{\mathbb{Q}}$. This is a contradiction. \endproof

\section{Exponential transcendence and exponentially transcendental powers}

In this section we recall some material from \cite{K} and then slightly modify a result from \cite{BKW}. 

\begin{defn}
A partial exponential field $(F,A,\exp)$ consists of a field $F$ of characteristic zero, a divisible subgroup $A$ of the additive group of $F$ and a homomorphism $\exp$ from $A$ to the multiplicative group of $F$.
\end{defn}

A total exponential field is a partial exponential field $(F,A,\exp)$ for which $F=A$. An $E$-ring is a commutative ring with identity, equipped with a homomorphism $E$ from its additive group to its multiplicative monoid. For variables $X_1,...,X_k$, $\mathbb{Z}[X_1,...,X_k]^E$ denotes the free $E$-ring generated by $\{X_1,...,X_k\}$. If $(F,A,\exp)$ is a partial exponential field then each $f\in \mathbb{Z}[X_1,...,X_k]^E$ determines a partial function $f:F^k\rightarrow F$ in the obvious way. The following definition, which comes from work of Macintyre \cite{Mac}, plays an important role in \cite{K}.
\begin{defn}
Let $(F,A,\exp)$ be a partial exponential field. Let $a\in F$ and $B\subseteq F$. We say that $a\in\text{ecl}(B)$ if there exist $f_1,...,f_k\in \mathbb{Z}[X_1,...,X_k,Y_1,...,Y_m]^E$, $a_1,...,a_k\in F$ and $b_1,...,b_m\in B$ such that 

\begin{enumerate}
\item $a=a_1$,

\item for all $i\in\{1,...,k\}$, $f_i(a_1,...,a_k,b_1,...,b_m)$ is defined and equal to zero and 

\item the determinant of the matrix $\left(\frac{\partial f_i}{\partial X_j}(a_1,...,a_k,b_1,...,b_m)\right)_{i,j\in\{1,...,k\}}$ is not zero.

\end{enumerate}

\end{defn}

The determinant of the Jacobian matrix will be defined at $(a_1,...,a_k,b_1,...,b_m)$ if each $f_i$ is and so there was no need to include that as an assumption. Kirby also considers another closure concept in \cite{K}. This is defined in terms of derivations which respect exponentiation. An exponential derivation $\partial$ on a partial exponential field $(F,A)$ is a map $\partial:F\rightarrow M$ such that $M$ is an $F$-vector space and the following conditions are satisfied.

\begin{enumerate}
\item For all $a,b\in F$, $\partial(a + b)=\partial(a) + \partial(b)$. 
\item For all $a,b\in F$, $\partial(a b)=a\partial(b) + b\partial(a)$. 
\item For all $a\in A$, $\partial(\exp(a))=\exp(a)\partial(a)$. 

\end{enumerate}

\begin{defn}
Let $(F,A,\exp)$ be a partial exponential field. Let $a\in F$ and $B\subseteq F$. We say that $a\in \text{ecl}'(B)$ if, for every exponential derivation $\partial$ on $(F,A,\exp)$, if $\partial(b)=0$ for all $b\in B$ then $\partial(a)=0$. 
\end{defn}

For the main results of \cite{K}, it is important to recall that a pregeometry $(X,\text{cl})$ consists of a set $X$ and a map $\text{cl}:\mathcal{P}(X)\rightarrow\mathcal{P}(X)$ such that the following conditions are satisfied. 

\begin{enumerate}

\item For all $B\subseteq X$, $B\subseteq \text{cl}(B)$ and $\text{cl}(\text{cl}(B))=\text{cl}(B)$. 

\item For all $B,C\subseteq X$, if $B\subseteq C$ then $\text{cl}(B)\subseteq \text{cl}(C)$.

\item For all $B\subseteq X$ and $a\in X$, if $a\in \text{cl}(B)$ then there is a finite subset $B_0\subseteq B$ such that $a\in \text{cl}(B_0)$. 

\item For all $B\subseteq X$ and $a,c\in X$, if $a\in \text{cl}(B\cup\{c\})\setminus\text{cl}(B)$ then $c\in \text{cl}(B\cup\{a\})$. 

\end{enumerate}

The following result is Theorem 1.1 in \cite{K}. It generalises earlier work of Wilkie \cite{W2}.

\begin{thm}\label{K1}
Let $(F,A,\exp)$ be a partial exponential field. For all $B\subseteq F$, $\text{\em ecl}(B)=\text{\em ecl}'(B)$. Furthermore, $(F,\text{\em ecl})$ is a pregeometry. 
\end{thm}

It is well-known that any pregeometry $(F,\text{cl})$ induces a notion of dimension. Let $a=(a_1,...,a_k)\in F^k$ and $B\subseteq F$. We say that $a$ is $\text{cl}$-independent over $B$ if, for all $i\in \{1,...,k\}$, $a_i\notin \text{cl}(B\cup\{a_1,...,a_{i-1},a_{i+1},...,a_k\})$. The dimension $\dim(a/B)$ of a finite tuple $a$ over a set $B$ is then defined to be the length of any maximal subtuple of $a$ which is $\text{cl}$-independent over $B$. We use some standard notational conventions regarding $\dim(a/B)$. For example, if $a=(a_1,...,a_k)$, $b=(b_1,...,b_m)$ and $c=(c_1,...,c_n)$ are finite tuples and $A$ and $B$ are sets then $\dim(a,b/A,B,c)=\dim(d/D)$, where $d=(a_1,...,a_k,b_1,...,b_m)$ and $D=A\cup B\cup\{c_1,...,c_n\}$. The following property of dimension is well-known and easy to prove.

\begin{fact}\label{fact}
Let $(F,\text{\em cl})$ be a pregeometry. Let $a\in F^k$, $c\in F^n$ and $B\subseteq F$. Then $$\dim(a,c/B)=\dim(c/B)+\dim(a/B,c).$$
\end{fact}

Given a field $F$ and a subfield $Q$, we have familiar pregeometries $(F,\text{acl})$ and $(F,\text{span}_Q)$, where $\text{acl}(B)$ is the algebraic closure in $F$ of the subfield generated by $B$ and $\text{span}_Q(B)$ is the $Q$-vector subspace of $F$ generated by $B$. The dimensions associated with these two pregeometries are denoted by $\text{td}$ and $\text{ldim}_Q$ respectively. Given a partial exponential field $(F,A)$, the dimension associated with $\text{ecl}$ is denoted by $\dim$. The notation $\dim(a)$, $\text{td}(a)$ and $\text{ldim}_Q(a)$ abbreviates $\dim(a/\emptyset)$, $\text{td}(a/\emptyset)$ and $\text{ldim}_Q(a/\emptyset)$.

Theorem 1.2 in \cite{K} applies to partial exponential fields, as well as total ones, provided all instances of exponentiation are defined, as can be seen from the proof. Making this explicit, the statement is as follows.

\begin{thm}\label{K2}
Let $(F,A,\exp)$ be a partial exponential field. Let $b\in A^N$ and $K\subseteq F$ be such that $\text{\em ecl}(K)=K$. Then $$\text{\em td}(b,\exp(b)/K)-\text{\em ldim}_{\mathbb{Q}}(b/K)\geq \dim(b/K).$$
\end{thm}

Note that $\exp(b)$ is an abbreviation of $(\exp(b_1),...,\exp(b_N))$, where $b=(b_1,...,b_N)$. The next result is a slight variant of Theorem 1.2 in \cite{BKW}. That theorem was stated for total exponential fields. The version we give here is suitable for our purposes and has almost the same proof. We work over an exponentially-algebraically closed parameter set $K\subseteq F$, which Bays, Kirby and Wilkie have in their proof though they omit it from the statement of the theorem.  

\begin{thm}\label{newBKW}
Let $(F,A,\exp)$ be a partial exponential field. Let $a=(a_1,...,a_N)\in A^N$ and $\lambda\in F$. Let $K\subseteq F$ be such that $K= \text{\em ecl}(K)$ and $\lambda\notin K$. Suppose $\lambda a=(\lambda a_1,...,\lambda a_N)\in A^N$ and assume $\exp(a)$ is multiplicatively independent. Then $$\text{\em td}(\text{\em exp}(a),\text{\em exp}(\lambda a)/K,\lambda)\geq N.$$
\end{thm}

\proof To prove this in the case where $F=A$, Bays, Kirby and Wilkie first show that $$\text{td}(\text{exp}(a,\lambda a)/K,\lambda)+\text{ldim}_{\mathbb{Q}(\lambda)}(a,\lambda a/\text{ker})-\text{ldim}_\mathbb{Q}(a,\lambda a/\text{ker})\geq 0$$ where $\text{ker}=\{c\in F:\exp(c)=1\}$. Their argument involves applying Theorem \ref{K2} (above) with $b=(a,\lambda a,\lambda)$ and so makes use of $\exp(\lambda)$. In our setting, $\exp(\lambda)$ might not be defined. We therefore modify the argument slightly so that Theorem \ref{K2} is applied with $b=(a,\lambda a)$. Adapting the first part of Section 4 in \cite{BKW}, we reason as follows.

Let $b=(a,\lambda a)$. By Theorem \ref{K2}, $$\text{td}(b,\exp(b)/K)-\text{ldim}_\mathbb{Q}(b/K)\geq 1$$ and so $$\text{td}(b/K)+\text{td}(\exp(b)/K,b)-\text{ldim}_\mathbb{Q}(b/K)\geq 1.$$ Now $$\text{td}(b/K)=\text{td}(b,\lambda/K)=\text{td}(\lambda/K)+\text{td}(b/K,\lambda)= 1+\text{td}(b/K,\lambda).$$ So then $$\text{td}(b/K,\lambda)+\text{td}(\exp(b)/K,b)-\text{ldim}_\mathbb{Q}(b/K)\geq 0.$$ We have $$\text{td}(b/K,\lambda)\leq \text{ldim}_{\mathbb{Q}(\lambda)}(b/K)$$ and $$\text{td}(\exp(b)/K,b)\leq\text{td}(\exp(b)/K,\lambda).$$ Therefore $$\text{td}(\exp(b)/K,\lambda)+\text{ldim}_{\mathbb{Q}(\lambda)}(b/K)-\text{ldim}_\mathbb{Q}(b/K)\geq 0.$$ Since $K=\text{acl}(K)$ and $\lambda\notin K$, we may deduce, exactly as in \cite{BKW}, that $$\text{ldim}_{\mathbb{Q}(\lambda)}(b/K)-\text{ldim}_\mathbb{Q}(b/K)\leq \text{ldim}_{\mathbb{Q}(\lambda)}(b/\text{ker})-\text{ldim}_\mathbb{Q}(b/\text{ker}).$$ This is done by applying part (iii) of Lemma 3.2 in \cite{BKW} and then Lemma 3.3 in \cite{BKW}. So we have $$\text{td}(\exp(b)/K,\lambda)+\text{ldim}_{\mathbb{Q}(\lambda)}(b/\text{ker})-\text{ldim}_\mathbb{Q}(b/\text{ker})\geq 0.$$ It remains to show that $$\text{ldim}_\mathbb{Q}(b/\text{ker})-\text{ldim}_{\mathbb{Q}(\lambda)}(b/\text{ker})\geq N.$$ This is done in Section 4 of \cite{BKW}, in the part subtitled ``Proof of Theorem 1.2", using no assumptions other than the ones we have made here. Note that $x$, $z$ and $n$, in the relevant part of \cite{BKW}, correspond to our $a$, $b$ and $N$, respectively. \endproof

\section{Proofs of the main results}

The first aim of this section is to prove Theorem \ref{newtheorem}. The full statement is as follows. 

\begin{thm}\label{strongertheorem}
Let $C_1,C_2\subseteq\mathbb{G}_m^N(\mathbb{C})$ be irreducible closed algebraic curves, with $N\geq 3$. Suppose there does not exist an algebraic subgroup $G\subseteq \mathbb{G}_m^N(\mathbb{C})$ of dimension $1$ such that $C_1\subseteq G$ and there does not exist an algebraic subgroup $H\subseteq\mathbb{G}_m^N(\mathbb{C})$ of dimension $2$ such that $C_1\cup C_2\subseteq H$. Further assume that at least one of the two curves is not defined over $\overline{\mathbb{Q}}$. Let $\mathcal{N}=\{n\in\mathbb{N}:[n]C_1\subseteq C_2\}$. Then $\bigcup\limits_{n\in \mathbb{N}\setminus\mathcal{N}}\{x\in C_1:x^n\in C_2\}$ is finite. 
\end{thm}

Let $\mathbb{R}_{(\mathbb{C},\exp)}$ denote the ordered field of real numbers expanded to include the usual exponential field structure on $\mathbb{C}=\mathbb{R}^2$. For some $\kappa>2^{\aleph_0}$, let $R_{(F,\exp)}$ be a $\kappa$-saturated and strongly $\kappa$-homogeneous elementary extension of $\mathbb{R}_{(\mathbb{C},\exp)}$, with $F=R^2$ and $\exp:F\rightarrow F$. Let $A\subseteq F$ be the subgroup of elements with infinitesimal modulus (including zero), so $A=\{z\in F:\left|z\right|<\frac{1}{n}\text{ for all }n\in\mathbb{N}\}$. Restricting $\exp$ to $A$ gives a partial exponential field $(F,A,\exp)$. We identify $\mathbb{C}$ with its image under the elementary embedding, so $\mathbb{C}\subseteq F$. Let $R_{(F,A,\exp,\mathbb{C})}$ be the expansion of $R_{(F,\exp)}$ obtained by adding unary predicates for $A$ and $\mathbb{C}$. Then let $\widetilde{R}_{(\widetilde{F},\widetilde{A},\exp,\widetilde{\mathbb{C}})}$ be the ultrapower of this with respect to a non-principal ultrafilter on $\mathbb{N}$. In particular, $(\widetilde{F},\widetilde{A},\exp)$ is an elementary extension of the partial exponential field $(F,A,\exp)$. If, in the above construction, we restrict $\exp$ to $D=\{z\in F:\left|z\right|\leq 1\}$, instead of to $A$, then neither $(F,D,\exp)$ nor $(\widetilde{F},\widetilde{D},\exp)$ is a partial exponential field, since $D$ is not a subgroup, but the structure $\widetilde{R}_{(\widetilde{F},\widetilde{D},\exp)}$ has the advantage of being o-minimal (see \cite{vdD}). 

For the rest of this section, $\text{acl}$ denotes algebraic closure in the field $\widetilde{F}$, $\text{ecl}$ denotes exponential-algebraic closure in the partial exponential field $(\widetilde{F},\widetilde{A},\exp)$ and $\text{dcl}$ denotes definable closure in the o-minimal structure $\widetilde{R}_{(\widetilde{F},\widetilde{D},\exp)}$.

\begin{lem}\label{two closures equal}
Let $B\subseteq \widetilde{\mathbb{C}}$. Then $\text{\em acl}(B)=\text{\em ecl}(B)$. 
\end{lem}

\proof Let $a\in \widetilde{F}$ and suppose $a\in \text{ecl}(B)$. Then there exist $f_1,...,f_k\in \mathbb{Z}[X_1,...,X_k,Y_1,...,Y_m]^E$, $a_1,...,a_k\in F$ and $b_1,...,b_m\in B$ such that $a=a_1$, $$f_1(a_1,...,a_k,b_1,...,b_m)=...=f_k(a_1,...,a_k,b_1,...,b_m)=0$$ and $(a_1,...,a_k)$ is an isolated solution of the system $$f_1(X_1,...,X_k,b_1,...,b_m)=...=f_k(X_1,...,X_k,b_1,...,b_m)=0.$$ Using the fact that $\widetilde{A}\cap \widetilde{\mathbb{C}}=\{0\}$, there are, for each $i\in\{1,...,k\}$, polynomials $f_{i,1},...,f_{i,k_i}\in \mathbb{Z}[X_1,...,X_k,Y_1,...,Y_m]$ such that the equation $f_i(X_1,...,X_k,b_1,...,b_m)=0$ is implied by the system $$f_{i,1}(X_1,...,X_k,b_1,...,b_m)=...=f_{i,k_i}(X_1,...,X_k,b_1,...,b_m)=0$$ and, for any $a_1',...,a_k'\in\widetilde{\mathbb{C}}$, $f_i(a_1',...,a_k',b_1,...,b_m)=0$ if and only if $$f_{i,1}(a_1',...,a_k',b_1,...,b_m)=...=f_{i,k_i}(a_1',...,a_k',b_1,...,b_m)=0.$$ It follows that if our tuple $(a_1,...,a_k)$ is in $\widetilde{\mathbb{C}}^k$ then $(a_1,...,a_k)$ is an isolated solution of $$\bigwedge\limits_{i=\{1,...,k\}}f_{i,1}(X_1,...,X_k,b_1,...,b_m)=...=f_{i,k_i}(X_1,...,X_k,b_1,...,b_m)=0$$ and so $a\in \text{acl}(B)$. 

Suppose $(a_1,...,a_k)\notin \widetilde{\mathbb{C}}^k$. Then, by basic o-minimality, $a_i\notin(\text{dcl}(B))^2$, for some $i\in\{1,...,k\}$. It follows that the complete type of $(a_1,...,a_k)$ over $B$, in the sense of $\widetilde{R}_{(\widetilde{F},\widetilde{D},\exp)}$, has infinitely many realisations. By o-minimality, this contradicts the fact that $(a_1,...,a_k)$ is an isolated solution of the system $$f_1(X_1,...,X_k,b_1,...,b_m)=...=f_k(X_1,...,X_k,b_1,...,b_m)=0.$$ \endproof

It is worth noting the following.

\begin{prop}\label{countable}
For each $B\subseteq \widetilde{R}$ and $a\in \widetilde{F}$, considered as a pair $a=(a_1,a_2)\in \widetilde{R}^2$, if $a\in \text{\em ecl}(B^2)$ then $a_1,a_2\in \text{\em dcl}(B)$.
\end{prop}

By Proposition \ref{relative implies}, to prove Theorem \ref{strongertheorem} it is sufficient to establish the following.

\begin{thm}\label{geometric}
Let $C_1,C_2\subseteq\mathbb{G}_m^3(\mathbb{C})$ be irreducible closed algebraic curves. Suppose there is no proper algebraic subgroup $H\subseteq \mathbb{G}_m^3(\mathbb{C})$ such that $C_1\subseteq H$. Further assume that at least one of the two curves is not defined over $\overline{\mathbb{Q}}$. Let $\mathcal{N}=\{n\in\mathbb{N}:[n]C_1\subseteq C_2\}$. Then $\bigcup\limits_{n\in \mathbb{N}\setminus\mathcal{N}}\{x\in C_1:x^n\in C_2\}$ is finite. 
\end{thm}

\proof Suppose $\bigcup\limits_{n\in \mathbb{N}\setminus\mathcal{N}}\{x\in C_1:x^n\in C_2\}$ is infinite. Let $e_1,...,e_k\in \mathbb{C}$ be algebraically independent such that $C_1$ and $C_2$ are defined over $K=\text{acl}(e_1,...,e_k)$. We may assume that, considered as a $2k$-tuple in the o-minimal structure $R_{(F,D,\exp)}$, $(e_1,...,e_k)$ is $\text{dcl}$-independent. Let $\sigma$ be an automorphism of $R_{(F,D,\exp)}$ such that $(\sigma(e_1),...,\sigma(e_k))$ is $\text{dcl}$-independent over $\mathbb{R}$ and $\sigma(e_i)-e_i\in A$, for all $i\in \{1,...,k\}$. Let $K'=\sigma(K)$, $C_1'=\sigma(C_1(F))$, $C_2'=\sigma(C_2(F))$ and $e_i'=\sigma(e_i)$, for all $i\in \{1,...,k\}$.

Letting $\widetilde{\mathbb{R}}$ denote the canonical extension of $\mathbb{R}$ to the ultrapower $\widetilde{R}_{(\widetilde{F},\widetilde{D},\exp)}$, we have that $(e'_1,...,e'_k)$, as a $2k$-tuple, is $\text{dcl}$-independent over $\widetilde{\mathbb{R}}$. Therefore, by Proposition \ref{countable}, the $k$-tuple $(e'_1,...,e'_k)$ is $\text{ecl}$-independent over $\widetilde{\mathbb{C}}$. Working in the underlying real-closed field of $R_{(F,D,\exp)}$, it is easy to see that, for each $n\in \mathbb{N}\setminus\mathcal{N}$ and $x\in C_1$, if $x^n\in C_2$ then there exists $z\in C_1'$ such that $z^n\in C_2'$, $x^{-1}z\in \exp(A)^3$ and $x^{-n}z^n\in \exp(A)^3$.

Letting $\widetilde{\mathbb{N}}$ denote the canonical extension of $\mathbb{N}$ to the ultrapower $\widetilde{R}_{(\widetilde{F},\widetilde{A},\exp,\widetilde{\mathbb{C}})}$, we have an infinite $n\in \widetilde{\mathbb{N}}\setminus\mathcal{N}$ and points $x\in C_1(\widetilde{\mathbb{C}})$ and $z\in C_1'(\widetilde{F})$ such that $x^n\in C_2(\widetilde{\mathbb{C}})$, $z^n\in C_2'(\widetilde{F})$, $x^{-1}z\in \exp(\widetilde{A})^3$, $x^{-n}z^n\in \exp(\widetilde{A})^3$ and $x$ is not algebraic over $K$. As elements of an ultrapower, $x$, $z$ and $n$ are really equivalence classes of sequences, so $x=[(x_j)_{j\in\mathbb{N}}]$, $z=[(z_j)_{j\in \mathbb{N}}]$ and $n=[(n_j)_{j\in\mathbb{N}}]$. By $x^n$ and $z^n$ we mean $[(x_j^{n_j})_{j\in\mathbb{N}}]$ and $[(z_j^{n_j})_{j\in\mathbb{N}}]$. 

The point $x^{-1}z$ lies on the curve $x^{-1}C_1'(\widetilde{F})$, while $x^{-n}z^n$ lies on $x^{-n}C_2'(\widetilde{F})$. These curves are defined over $K'(x,x^n)$. We consider various cases and obtain a contradiction in each.

Case 1: Suppose $x^{-1}z$ is multiplicatively independent. Let $u$ be a maximal subtuple of the 6-tuple $(x,x^n)$ such that $u$ is $\text{acl}$-independent over $K\cup\{n\}$. It follows by Theorem \ref{newBKW} that $n\in \text{ecl}(K,K',u)$. Then $n\in \text{ecl}(K,u)$, since $e_1',...,e_k'$ are $\text{ecl}$-independent over $\widetilde{\mathbb{C}}$. By Proposition \ref{two closures equal}, $n\in \text{acl}(K,u)$. Since $n$ is infinite, $n\notin \text{acl}(K)$. This contradicts the $\text{acl}$-independence of $u$ over $K\cup\{n\}$. 

Case 2: Suppose $x^{-1}z$ has multiplicative dimension $2$. Then there is an algebraic subgroup $H\subseteq\mathbb{G}_m^3(\widetilde{F})$ of dimension $2$ such that $x^{-1}z\in H$ and $x^{-n}z^n\in H$. Then $x\in zH$ and $x^n\in z^nH$. 

Subcase 2a: Suppose both $C_1(\widetilde{F})\subseteq zH$ and $C_2(\widetilde{F})\subseteq z^nH$. This contradicts the fact that $\bigcup\limits_{m\in\mathbb{N}\setminus\mathcal{N}}\{x\in C_1:x^m\in C_2\}$ is infinite, unless $zH$ is a translate of $H$ by a torsion point, which would contradict the fact that $C_1$ is not contained in a proper algebraic subgroup of $\mathbb{G}_m^3(\mathbb{C})$.

Subcase 2b: Suppose $C_2(\widetilde{F})\nsubseteq z^nH$. Then $x^n\in \text{acl}(z^n,K,K')$. It follows from Corollary \ref{Liardet} and our assumptions that $x^n\notin \text{acl}(K)$. Since $e'_1,...,e'_k$ are $\text{acl}$-independent over $\widetilde{\mathbb{C}}$, $x^n\notin\text{acl}(K,K')$. Therefore, by exchange, $z^n\in \text{acl}(K,K',x^n)$. As in Case 1, let $u$ be an $\text{acl}$-basis for $(x,x^n)$ over $K\cup\{n\}$. By Theorem \ref{newBKW} and the fact that $\text{td}(x^{-1}z,x^{-n}z^n/K,K',u,n)\leq 1$, $n\in \text{ecl}(K,K',u)$. We obtain a contradiction as in Case 1.

Subcase 2c: The subcase where $C_1(\widetilde{F})\nsubseteq zH$ is analogous to the previous one. This time we use the fact that $x\notin \text{acl}(K)$. 

Case 3: Suppose $x^{-1}z$ has multiplicative dimension $1$. Then there is an algebraic subgroup $H\subseteq \mathbb{G}_m^3(\widetilde{F})$ of dimension $1$ such that $x\in zH$ and $x^n\in z^nH$. By Proposition \ref{C not a coset}, $C_1(\widetilde{F})\nsubseteq zH$ and $C_2(\widetilde{F})\nsubseteq z^nH$. So $x\in \text{acl}(K,K',z)$ and $x^n\in \text{acl}(K,K',z^n)$. By reasoning from Subcases 2b and 2c, we then have $z\in \text{acl}(K,K',x)$ and $z^n\in \text{acl}(K,K',x^n)$. Let $u$ be an $\text{acl}$-basis for $(x,x^n)$ over $K\cup\{n\}$. By Theorem \ref{newBKW} and the fact that $\text{td}(x^{-1}z,x^{-n}z^n/K,K',u,n)=0$, $n\in \text{ecl}(K,K',u)$. We obtain a contradiction as in Case 1.

Case 4: Suppose $x^{-1}z$ has multiplicative dimension zero. Then it is the point $(1,1,1)$ and the curves $C_1(F)$ and $C_1'(F)$ have infinitely many points in common. It follows that $C_1(F)=C_1'(F)$. Similarly, $C_2(F)=C_2'(F)$. Then $C_1(F)$ and $C_2(F)$ are both defined over $K\cap K'=\overline{\mathbb{Q}}$, which is a contradiction. \endproof

We now prove Theorem \ref{higher dim}, stated again here for convenience.

\begin{thm}\label{higher dim repeat}
Let $V_1, V_2\subseteq\mathbb{G}_m^N(\mathbb{C})$ be irreducible closed algebraic subvarieties, with $\dim(V_1)\geq 1$ and $\dim(V_1)+\dim(V_2)< N$. Suppose that, for any $\mathbb{Q}$-linearly independent tuples of integers $$(n_{1,1},...,n_{1,N}),...,(n_{k,1},...,n_{k,N}),$$ with $k=\dim(V_1)+1$, and any generic point $(x_1,...,x_N)$ of $V_1$, the monomials $$x_1^{n_{1,1}}\cdot...\cdot x_N^{n_{1,N}},...,x_1^{n_{k,1}}\cdot...\cdot x_N^{n_{k,N}}$$ are algebraically independent.  Let $\mathcal{N}=\{n\in\mathbb{N}:[n]V_1\subseteq V_2\}$. Then $\bigcup\limits_{n\in\mathbb{N}\setminus\mathcal{N}}\{x\in V_1:x^n\in V_2\}$ is not Zariski-dense in $V_1$. 
\end{thm}

\proof Suppose $\bigcup\limits_{n\in\mathbb{N}\setminus\mathcal{N}}\{x\in V_1:x^n\in V_2\}$ is Zariski-dense in $V_1$. Let $e_1,...,e_k\in \mathbb{C}$ and $e'_1,...,e'_k\in F$ be as in the previous proof, except with $V_1$ in place of $C_1$ and $V_2$ in place of $C_2$. Let $K=\text{acl}(e_1,...,e_k)$ and $K'=\text{acl}(e'_1,...,e'_k)$. We obtain $V'_1$ and $V'_2$ analogously to how $C'_1$ and $C'_2$ were obtained. 

Let $(x_j)_{j\in\mathbb{N}}$ and $(n_j)_{j\in\mathbb{N}}$ be sequences in $V_1$ and $\mathbb{N}\setminus\mathcal{N}$ respectively such that $x_j^{n_j}\in V_2$, for all $j\in \mathbb{N}$, and $\{x_j:j\in \mathbb{N}\}$ is Zariski-dense in $V_1$. Then $(n_j)_{j\in\mathbb{N}}$ is not eventually constant. We get a generic (over $K$) point $x=[(x_j)_{j\in\mathbb{N}}]\in V_1(\widetilde{\mathbb{C}})$, an infinite $n=[(n_j)_{j\in\mathbb{N}}]\in \widetilde{\mathbb{N}}\setminus\widetilde{\mathcal{N}}$ and some $z=[(z_j)_{j\in\mathbb{N}}]\in V_1'(\widetilde{F})$ such that $x^n\in V_2(\widetilde{\mathbb{C}})$, $z^n\in V_2'(\widetilde{F})$, $x^{-1}z\in \exp(\widetilde{A})^N$ and $x^{-n}z^n\in \exp(\widetilde{A})^N$. We may assume $z$ is a generic point of $V'_1(\widetilde{F})$ over $K'$ and that $\text{td}(z_j/K')=\text{td}(z_j/K,K',x_j)$ for all $j\in\mathbb{N}$. 

Suppose $x^{-1}z$ satisfies $d$ independent multiplicative relations, with $d\leq \dim(V_1)+1$. Then there is a surjective algebraic homomorphism $\phi:\mathbb{G}_m^N(\widetilde{F})\rightarrow\mathbb{G}_m^d(\widetilde{F})$ such that $\phi(x)^{-1}=\phi(z)$. It follows that $\phi(x),\phi(z)\in \mathbb{G}_m^d(\widetilde{L})$, where $L=\overline{\mathbb{Q}}$. Suppose $\phi^{-1}(\phi(z))\cap V_1'(\widetilde{F})$ has dimension greater than $\dim(V_1')-d$. Then the Zariski-closure of $\phi(V_1'(\widetilde{F}))$ has dimension less than $d$ and, since $\mathbb{G}_m^d(\overline{\mathbb{Q}})$ is Zariski-dense in it, is defined over $\overline{\mathbb{Q}}$. This contradicts our assumptions. So $$0\leq\dim(\phi^{-1}(\phi(z))\cap V_1'(\widetilde{F}))\leq\dim(V_1')-d=\dim(V_1)-d.$$ Therefore $d\leq\dim(V_1)$. So $x^{-1}z$ satisfies at most $\dim(V_1)$ independent multiplicative relations. We may assume $d$ is maximal.

Let $u$ be a maximal subtuple of $(x,x^n)$ such that $u$ is $\text{acl}$-independent over $K\cup\{n\}$. Then $\text{td}(x^{-1}z,x^{-n}z^n/K,K',u,n)\leq \dim(V_1)-d+\dim(V_2)$ and the multiplicative dimension of $x^{-1}z$ is at least $N-d$. By Theorem \ref{newBKW}, $n\in \text{ecl}(K,K',u)$ and so $n\in \text{ecl}(K,u)=\text{acl}(K,u)$.  Since $n$ is infinite, $n\notin \text{acl}(K)$. This contradicts the $\text{acl}$-independence of $u$ over $K\cup\{n\}$. \endproof

It follows easily from known results that the set $\mathcal{N}$ in this theorem is finite.

\bibliographystyle{abbrv}

\end{document}